\documentclass[11pt]{article}
\usepackage{graphicx}
\usepackage{amsmath}
\usepackage{epsfig}  
\usepackage{verbatim}  
\usepackage{tabularx}
\usepackage{epic, eepic}
\usepackage{xcolor}
\usepackage{lineno}
\usepackage{amssymb,latexsym,amscd,amsmath,amsfonts,enumerate,supertabular}
\usepackage{subfigure}
\usepackage{longtable}
\newtheorem{Theorem}{Theorem}[section] 
\newtheorem{Proposition}[Theorem]{Proposition} 

\newtheorem{remark}[Theorem]{Remark}

\begin{document}

\title{{\bf A unified approach to study the existence and numerical solution of functional differential equation }}

\author{ Dang Quang A$^{\text a}$,  Dang Quang Long$^{\text b}$\\
$^{\text a}$ {\it\small Center for Informatics and Computing, VAST}\\
{\it\small 18 Hoang Quoc Viet, Cau Giay, Hanoi, Vietnam}\\
{\small Email: dangquanga@cic.vast.vn}\\
$^{\text b}$ {\it\small Institute of Information Technology, VAST,}\\
{\it\small 18 Hoang Quoc Viet, Cau Giay, Hanoi, Vietnam}\\
{\small Email: dqlong88@gmail.com}}
\date{ }
%
\date{}          
\maketitle
\begin{abstract}
In this paper we consider a class of boundary value problems for third order nonlinear functional differential equation.  By the reduction of the problem to operator equation we establish the existence and uniqueness of solution and construct a numerical method for solving it. We prove that the method is of second order accuracy and obtain an estimate for total error.  Some examples demonstrate the validity of the obtained theoretical results and the efficiency of the numerical method. The approach used for the third order nonlinear functional differential equation
can be applied to functional differential equations of any orders.

\end{abstract}
{\small
}
\section{Introduction}\label{sec1}

Functional differential equations  have numerous applications in engineering and sciences \cite{Hale}. Therefore, for the last decades they have been studied by many authors. There are many works concerning the numerical solution of both initial and boundary value problems for them. The methods used are diverse including collocation method \cite{Reut}, iterative methods \cite{Bica,Khuri}, neural networks \cite{Hou,Raja}, and so on. Below we mention some results of typical works.\par 
First it is  worthy to mention the work of Reutskiy in 2015 \cite{Reut}. In this work the author considered the linear pantograph  functional differential equation with proportional delay
\begin{align*}
u^{(n)}=\sum_{j=0}^{J}\sum _{k=0}^{n-1} p^{jk}(x) u^{(k)} (\alpha _{j}x)+f(x), \quad x\in [0,T]
\end{align*}
associated with initial or boundary conditions. Here $\alpha _j$ are constants ($0< \alpha _j <1$). The author proposed a  method, where  the initial equation is replaced by  an approximate equation which has an exact analytic solution with a set of free parameters. These free parameters are determined by the use of the collocation procedure. Many examples show the efficiency of the method but no errors estimates are obtained.\par

In 2016 Bica et al. \cite{Bica} considered the boundary value problem (BVP)
\begin{equation}\label{eq:1}
\begin{split}
x^{(2p)}(t)=f(t,x(t), x(\varphi (t))), \quad t\in [a,b],\\
x^{(i)}(a)=a_i, \; x^{(i)}(b)=b_i,\quad i=\overline{0, p-1}
\end{split}
\end{equation}
where $\varphi : [a,b] \rightarrow \mathbb{R}, \; a \le \varphi (t) \le b, \forall t\in [a,b]$. For solving the problem, the authors constructed successive approximations for the  equivalent integral equation with the use of cubic spline interpolation at each iterative step. The error estimate was obtained for the approximate solution under the very strong conditions including $(\alpha +13 \beta )(b-a) M_G <1$, where $\alpha$ and  $\beta $ are the Lipshitz coefficients of the function $f(s,u,v)$ in the variables $u$ and $v$, respectively; $M_G$ is a number such that $|G(t,s)| \le M_G \; \forall t,s \in [a,b]$, $G(t,s)$ being the Green function for the above problem. Some numerical experiments demonstrate the convergence of the  proposed iterative method. But is a regret that in the proof of the error estimate for fourth order nonlinear BVP there is a vital mistake when the authors by default considered that the partial derivatives $\frac{\partial ^3 G}{\partial s^3}, \frac{\partial ^4 G}{\partial s^4}$ are continuous in $[a,b] \times [a,b]$. But it is invalid because $\frac{\partial ^3 G}{\partial s^3}$ has discontinuity on the line $s=t$. Due to this mistake the authors obtained that the error of the method for fourth order BVP is $O(h^4)$. Although in \cite{Bica} the method is constructed for the general function $\varphi (t)$ but in all numerical examples only the particular case $\varphi (t) =\alpha  t$ was considered and the conditions of convergence were not verified.\par 
Recently, in 2018 Khuri and Sayfy \cite{Khuri} proposed a Green function based iterative method for functional differential equations of arbitrary orders. But the scope of application of the method is very limited due to the difficulty in calculation of integrals at each iteration.\par
For  solving functional differential equations, besides analytical and  numerical methods recently computational intelligence algorithms also are used (see, e.g.,\cite{Hou,Raja}), where feed-forward artificial neural networks of different architecture are applied. These algorithms are heuristic, so no errors estimates are obtained and they require large computational efforts.\par 
In this paper we propose a new approach to functional differential equations (FDE). Although this approach can be applied to functional differential equations of any orders with nonlinear terms containing derivatives but for simplicity we consider the FDE of the form
\begin{equation}\label{eq:2}
u^{\prime \prime \prime}=f(t,u(t), u(\varphi (t))), \quad t\in [0,a]
\end{equation}
associated with the general boundary conditions
\begin{equation}\label{eq:3}
 \begin{split}
 B_1[u]=\alpha_1 u(0)+\beta_1 u'(0) + \gamma_1 u''(0) =b_1,\\
 B_2[u]=\alpha_2 u(0)+\beta_2 u'(0) + \gamma_2 u''(0) =b_2,\\
 B_3[u]=\alpha_3 u(1)+\beta_3 u'(1) + \gamma_3 u''(1) =b_3,\\
 \end{split}
 \end{equation}
or
 \begin{equation}\label{eq:4}
 \begin{split}
 B_1[u]=\alpha_1 u(0)+\beta_1 u'(0) + \gamma_1 u''(0) =b_1,\\
 B_2[u]=\alpha_2 u(1)+\beta_2 u'(1) + \gamma_2 u''(1) =b_2,\\
 B_3[u]=\alpha_3 u(1)+\beta_3 u'(1) + \gamma_3 u''(1) =b_3,\\
 \end{split}
 \end{equation}
 such that

 \begin{equation*}
 Rank
\begin{pmatrix}
\alpha_1 & \beta_1 & \gamma_1 & 0 &0&0\\
\alpha_2 & \beta_2 & \gamma_2 &0 &0&0\\
0& 0 & 0 & \alpha_3& \beta_3 & \gamma_3\\
\end{pmatrix}=3.
\end{equation*}
As in \eqref{eq:1}, the function $\varphi (t)$ is assumed to be continuous and maps $[0,a]$ into itself. \par Developing the unified approach for fully third order nonlinear differential equation 
\begin{align*}
u^{\prime \prime \prime}=f(t,u(t),u^{\prime}(t),u^{\prime \prime}(t)
\end{align*}
in the previous works \cite{Dang1,Dang2}, in this paper we establish the existence and uniqueness of solution of the problem \eqref{eq:2}-\eqref{eq:3}  and propose an iterative method for finding the solution at both continuous and discrete levels. Some examples demonstrate the validity of obtained theoretical results and the efficiency of the proposed numerical method.

\section{Existence and uniqueness of solution}\label{sec2}
Following the approach in \cite{Dang1,Dang2} (see also \cite{Dang3,Dang4}) to investigate the problem \eqref{eq:2}-\eqref{eq:3} we introduce the nonlinear operator $A$ defined in the space of continuous functions $C[0,a]$ by the formula:
\begin{equation}\label{eq:5}
(A\psi )(t) =f(t,u(t), u(\varphi (t))),
\end{equation}
where $u(t)$ is the solution of the problem
\begin{equation}\label{eq:6}
\begin{split}
u'''(t)&=\psi (t), \quad 0<t<1\\
B_1[u]&=b_1, B_2[u]=b_2, B_3[u]=b_3,
\end{split}
\end{equation}
where $B_1[u], B_2[u], B_3[u]$ are defined by \eqref{eq:3}. It is easy to verify the following
\begin{Proposition}\label{Prop1}
If the function $\psi$ is a fixed point of the operator $A$, i.e., $\psi$ is the solution of the operator equation 
\begin{equation}\label{eq:7}
A\psi = \psi ,
\end{equation} 
where $A$ is defined by \eqref{eq:5}-\eqref{eq:6}
then the function $u(t)$ determined from the BVP \eqref{eq:6} is a solution of the BVP \eqref{eq:2}-\eqref{eq:3}. Conversely, if the function $u(x)$ is the solution of the BVP \eqref{eq:2}-\eqref{eq:3} then the function 
\begin{equation*}
\psi (t)= f(t,u(t), u(\varphi (t)) )
\end{equation*}
satisfies the operator equation \eqref{eq:7}.
\end{Proposition}
Now, let $G(t,s) $ be the Green function of the problem \eqref{eq:6}. Then the solution of the problem can be represented in the form
\begin{equation}\label{eq:8}
u(t)= g(t)+ \int_0^a G(t,s)\psi (s) ds,
\end{equation}
where $g(t) $ is the polynomial of second degree satisfying the boundary conditions
\begin{equation}\label{eq:9}
B_1[g]=b_1, B_2[g]=b_2, B_3[g]=b_3,
\end{equation}
Denote 
\begin{equation}\label{eq:10}
M_0  =\max_{0\leq t\leq a} \int_0^1 |G(t,s)|ds .
\end{equation}
For any positive number $M$ define the domain
\begin{equation}\label{eq:11}
\mathcal{D}_M= \Big \{ (t,u,v) \mid  0\leq t\leq a; |u|\leq \|g\|+M_0 M; |v|\leq \|g\|+M_0 M \Big \},
\end{equation}
where $\|g\|= \max _{0\le t\le a}|g(t)|$.\\
As usual, we denote by $B[0,M]$ the closed ball of the radius $M$ centered at $0$ in the space of continuous functions $C[0,a]$. 
\begin{Theorem}\label{thm1}
Assume that:
\begin{description}
\item (i) The function $\varphi (t)$ is a continuous map from $[0,a]$ to $[0,a]$.
\item (ii) The function $f(t,u,v)$ is continuous and bounded by $M$ in the domain $\mathcal{D_M}$, i.e.,
\begin{equation}\label{eq:12}
|f(t,u,v)|\le M \quad \forall (t,u,v) \in \mathcal{D_M}.
\end{equation}
\item (iii) The function $f(t,u,v)$ satisfies the Lipshitz conditions in the variables $u,v$ with the coefficients $L_1, L_2 \ge 0$ in $\mathcal{D}_M$, i.e.,
\begin{equation}\label{eq:13}
\begin{split}
|f(t,u_2,v_2)-f(t,u_1,v_1)|\le L_1 |u_2-u_1|+L_2 |v_2-v_1| \quad \\
\forall (t,u_i,v_i)\in \mathcal{D}_M\;
(i=1,2)
\end{split}
\end{equation}
\item (iv)   
\begin{equation}\label{eq:q}
q:= (L_1+L_2)M_0 <1.
\end{equation}

\end{description}
The the problem \eqref{eq:2}-\eqref{eq:3} has a unique solution $u(t) \in C^3[0,a]$, satisfying the estimate
\begin{equation}\label{eq:14}
|u(t)| \le \|g\| +M_0M \quad \forall t\in[0,a].
\end{equation}
\end{Theorem}
\noindent {\bf Proof.} The proof of the theorem will be done in the following steps:\\
First we show that the operator $A$ is a mapping  $B[0,M]\rightarrow B[0,M]$. Indeed, for any $\psi \in B[0,M]$ we have $\|\psi \| \le M$. Let $u(t)$ be the solution of the problem \eqref{eq:6}. From \eqref{eq:8} it follows 
\begin{equation}\label{15}
|u(t)| \le \|g\| +M_0M \quad \forall t\in [0,a].
\end{equation}
Since $0 \le \varphi (t) \le a$ we also have
\begin{equation*}
|u(\varphi (t)| \le \|g\| +M_0M \quad \forall t\in [0,a].
\end{equation*}
Therefore, if $t\in [0,a]$ then $(t,u(t), u(\varphi (t)))\in \mathcal{D}_M$. By the assumption \eqref{eq:12} we have
$|f(t,u(t), u(\varphi (t)))| \le M \quad \forall t\in [0,a]$. In view of \eqref{eq:5} we have $|(A\psi)(t)|\le M \quad \forall t\in [0,a]$. It means $\| A\psi\| \le M$ or $A\psi \in B[0,M]$.\\
Next, we prove that $A$ is a contraction in $B[0,M]$. Let $\psi _1, \psi_2 \in B[0,M]$ and $u_1(t), u_2(t)$ be the solutions of the problem \eqref{eq:6}, respectively. Then from the assumption \eqref{eq:13} we obtain
\begin{equation}\label{eq:16}
|A\psi_2 -A\psi_1| \le L_1 |u_2(t)-u_1(t)| +L_2|u_2(\varphi (t))-u_1(\varphi (t))|.
\end{equation}
From the representations
\begin{equation*}
u_i(t)= g(t)+ \int_0^a G(t,s)\psi_i (s) ds,\quad (i=1,2)
\end{equation*}
and \eqref{eq:10}
it is easy to obtain
\begin{align*}
|u_2(t)-u_1(t)|& \le M_0 \|\psi_2 -\psi_1\|,\\
|u_2(\varphi  (t))-u_1(\varphi (  t))|& \le M_0 \|\psi_2 -\psi_1\|
\end{align*}
Combining the above estimates and \eqref{eq:16}, in view of the assumption \eqref{eq:q} we obtain
\begin{align*}
\|A\psi _2 -A\psi_1\|\le q\| \psi_2 -\psi_1\|, \quad q<1.
\end{align*} 
Thus, $A$ is a contraction mapping in $B[0,M]$.\par
Therefore, the operator equation \eqref{eq:7} has a unique solution $\psi \in B[0,M]$. By Proposition \ref{Prop1} the solution of the problem \eqref{eq:6} for this right-hand side $\psi (t)$ is the solution of the original problem \eqref{eq:2}-\eqref{eq:3}.
\section{Solution method and its convergence}\label{sec3}
Consider the following iterative method:
\begin{enumerate}
\item Given $\psi_0 \in B[0,M]$, for example,
\begin{equation}\label{iter1c}
\psi_0(t)=f(t,0,0).
\end{equation}
\item Knowing $\psi_k(t)$  $(k=0,1,...)$ compute
\begin{equation}\label{iter2c}
\begin{split}
u_k(t) &= g(t)+\int_0^a G(t,s)\psi_k(s)ds ,\\
v_k(t) &= g(\varphi (t))+\int_0^1 G(\varphi(t),s)\psi_k(s)ds ,
\end{split}
\end{equation}
\item Update
\begin{equation}\label{iter3c}
\psi_{k+1}(t) = f(t,u_k(t),v_k(t)).
\end{equation}
\end{enumerate}

Set 

\begin{equation}\label{eqpkd}
p_k=\dfrac{q^k}{1-q} ,\; d=\|\psi _1 -\psi _0\|.
\end{equation} 
\begin{Theorem}[Convergence]\label{thm2} Under the assumptions of Theorem \ref{thm1} the above iterative method converges and there holds the estimate
\begin{equation*}
\|u_k-u\| \leq M_0p_kd, 
\end{equation*} 
where $u$ is the exact solution of the problem \eqref{eq:2}-\eqref{eq:3} and $M_0$ is given by \eqref{eq:10}.
\end{Theorem}
This theorem follows straightforward from the convergence of the successive approximation method for finding fixed point of the operator $A$ and the representations \eqref{eq:8} and the first equation in \eqref{iter2c}.\par

To numerically realize the above iterative method we construct the corresponding discrete iterative method. For this purpose cover the interval $[0, a]$   by the uniform grid $\bar{\omega}_h=\{t_i=ih, \; h=a/N, i=0,1,...,N  \}$ and denote by $\Phi_k(t), U_k(t), V_k(t)$ the grid functions, which are defined on the grid $\bar{\omega}_h$ and approximate the functions $\psi_k (t), u_k(t), v_k(t)$ on this grid, respectively.\par
Below we describe the discrete iterative method:\\
\begin{enumerate}
\item Given 
\begin{equation}\label{iter1d}
\Psi_0(t_i)=f(t_i,0,0),\ i=0,...,N. 
\end{equation}
\item Knowing $\Psi_k(t_i),\; k=0,1,...; \; i=0,...,N, $  compute approximately the definite integrals \eqref{iter2c} by the trapezoidal rule
\begin{equation}\label{iter2d}
\begin{split}
U_k(t_i) &=g(t_i) +\sum _{j=0}^N h\rho_j G(t_i,t_j)\Psi_k(t_j), \\
V_k(t_i) &= g(\xi_i)+\sum _{j=0}^N h\rho_j G(\xi_i,t_j)\Psi_k(t_j) ,  i=0,...,N,
\end{split}
\end{equation}
\noindent where $\rho_j$ are the weights 
\begin{equation*}
\rho_j = 
\begin{cases}
1/2,\; j=0,N\\
1, \; j=1,2,...,N-1
\end{cases}
\end{equation*}
and $\xi_i=\varphi (t_i)$.
\item Update
\begin{equation}\label{iter3d}
\Psi_{k+1}(t_i) = f(t_i,U_k(t_i),V_k(t_i)).
\end{equation}
\end{enumerate}
Now study the convergence of the above discrete iterative method. For this purpose we need some auxiliary results.
\begin{Proposition}\label{prop2}
If the function $f(t,u,v)$ has all partial derivatives continuous up to second order and the function $\varphi (t)$ also has continuous derivatives  up to second order then the functions $\psi _k(t), u_k(t), v_k(t)$ constructed by the iterative method \eqref{iter1c}-\eqref{iter3c} also have continuous derivatives  up to second order.
\end{Proposition}
This proposition is obvious.
\begin{Proposition}\label{prop3}
For any function $\psi (t) \in C^2[0,a]$ there hold the estimates
\begin{align}
\int_0^a G (t_i,s) \psi (s) ds = \sum _{j=0}^N h\rho_j G(t_i,s_j)\psi(s_j) +O(h^2), \label{eq:p3.1}\\ 
\int_0^a G (\xi_i,s) \psi (s) ds = \sum _{j=0}^N h\rho_j G(\xi_i,s_j)\psi(s_j) +O(h^2), \label{eq:p3.2}
\end{align}
where in order to avoid possible confusion we denote $s_j=t_j$.
\end{Proposition}
\noindent{\it Proof.} The validity of \eqref{eq:p3.1} is guaranteed by \cite[Proposition 3]{Dang2}. Here we notice that \eqref{eq:p3.1} is not automatically deduced from the estimate for the composite trapezoidal rule because the function $\frac{\partial ^2G(t_i,s)}{\partial s^2}$ has discontinuity at $s=t_i$.\par 
Now we prove the estimate \eqref{eq:p3.2}. Since $0\le \xi_i =\varphi (t_i) \le a$, there are two cases.\\
Case 1: $\xi_i$ coincides with one node $s_j$ of the grid $\bar{\omega_h}$, i.e., there exists $s_j \in \bar{\omega_h}$ such that $\xi_i=s_j$. Because the Green function $G(t,s)$ as a function of $s$, it is continuous function at $s=\xi_i$ and is a polynomial of $s$ in the intervals $[0,\xi_i$ and $[\xi_i,a]$, we have
\begin{equation*}
\begin{aligned}
&\int_0 ^a G(\xi_i ,s)\psi (s) ds = \int_0^{\xi_i} G(\xi_i ,s)\psi (s) ds +\int_{\xi_i}^a G(\xi_i ,s)\psi (s) ds \\
&= h\Big (  \frac{1}{2}G(\xi_i,s_0)\psi(s_0) + \sum_{m=1}^{j-1} G(\xi_i,s_m) \psi(s_m)+ \frac{1}{2}G(\xi_i,s_j)\psi(s_j) \Big ) +O(h^2)\\
& +h\Big (  \frac{1}{2}G(\xi_i,s_j)\psi(s_j) + \sum_{m=j+1}^{N-1} G(\xi_i,s_m) \psi(s_m)+ \frac{1}{2}G(\xi_i,s_N)\psi(s_N) \Big ) +O(h^2)\\
&=\sum _{j=0}^N h\rho_j G(t_i,s_j)\psi(s_j) +O(h^2).
\end{aligned}
\end{equation*}
Thus, \eqref{eq:p3.2} is proved for Case 1.\\
Case 2: $\xi_i$ lies between $s_l$ and $s_{l+1}$, i.e., $s_l <\xi_i <s_{l+1}$ for some $l =\overline{0, N-1}$. In this case, we represent
\begin{equation}\label{eq:27}
\int_0 ^a G(\xi_i ,s)\psi (s) ds = \int _0 ^{s_l} F(s)ds +\int _{s_l} ^{\xi_i} F(s)ds+
\int _{\xi_i} ^{s_{l+1}} F(s)ds +\int_{s_{l+1}} ^{a} F(s)ds.
\end{equation}
Here, for short we denote $F(s)=G(\xi_i,s)\psi (s)$. Note that $F(s) \in C^2 $ in $[s_{l}, \xi_i]$ and $[\xi_i, s_{l+1}]$. Applying the composite trapezoidal rule to the first and the last integrals in the right-hand side of  \eqref{eq:27}  we obtain
\begin{equation}\label{eq:28}
\begin{aligned}
T_1:&=  \int _0 ^{s_l} F(s)ds +\int_{s_{l+1}} ^{a} F(s)ds\\
&= \sum _{j=0}^{l}\rho _{j}^{(l-)}F(s_{j}) + \sum _{j=l+1}^N \rho _{j}^{(l+)}F(s_{j})+O(h^2),
\end{aligned}
\end{equation}
where
\begin{align*}
\rho _{j}^{(l-)} =  \left \{  \begin{array}{ll}
\frac{1}{2}, \quad j=0,l\\
1, \quad, 1<j<l
 \end{array} \right. , \quad 
 \rho _{j}^{(l+)} =  \left \{  \begin{array}{ll}
\frac{1}{2}, \quad j=l+1,N\\
1, \quad l+1<j<N
 \end{array} \right.
\end{align*}
For calculating the second and the third integrals in the right-hand side of \eqref{eq:27} we use the trapezoidal rule

\begin{equation}\label{eq:28a}
\begin{aligned}
T_2:&=  \int _{s_l} ^{\xi_i} F(s)ds+ \int _{\xi_i} ^{s_{l+1}} F(s)ds \\
&= \frac{1}{2} \big [ (F(s_l)+F(\xi_i))(\xi_i-s_l) + (F(\xi_i)+F(s_{l+1}))(s_{l+1}-\xi_i)\big ]  +O(h^2).
\end{aligned}
\end{equation}
Using the points $s_l$ and $s_{l+1}$ for linearly interpolating $F(s)$ in the point $\xi_i$ we have
\begin{align*}
F(\xi_i)=F(s_l)\frac{\xi_i -s_{l+1}}{s_l-s_{l+1}}+F(s_{l+1})\frac{\xi_i-s_l}{s_{l+1}-s_l}+O(h^2).
\end{align*}
From here we obtain
\begin{equation}\label{eq:29}
F(\xi_i) (s_{l+1}-s_l)=F(s_l)(s_{l+1}-\xi_i)+F(s_{l+1})(\xi_i-s_l)+O(h^3).
\end{equation}
Now, transforming $T_2$ we have
\begin{align*}
T_2= \frac{1}{2} \big [ F(s_l) (\xi_i-s_l) + F(s_{l+1}) (s_{l+1}-\xi_i)  \big ]
 +F(\xi_i)(s_{l+1}-s_l)+O(h^2)
\end{align*}
Further, in view of \eqref{eq:29} it is easy to obtain
\begin{equation*}
T_2=\frac{1}{2}h(F(s_l)+F(s_{l+1}))+ O(h^3).
\end{equation*}
Taking into account the above estimate, \eqref{eq:28} and \eqref{eq:27} we have
\begin{equation*}
\int _0^a G(\xi_i,s)\psi(s) ds =\sum_{j=0}^N h\rho_j G(\xi_i,s_j)\psi(s_j)+O(h^2).
\end{equation*}
Thus, \eqref{eq:p3.2} is proved for Case 2 and the proof of Proposition \ref{prop3} is complete.

\begin{remark}
If in Proposition \ref{prop3} replace $G(t_i,s)$ and $G(\xi_i, s)$ by $|G(t_i,s)|$ and $|G(\xi_i, s)|$, respectively then we obtain the analogous estimates
\begin{align}
\int_0^a |G (t_i,s)| \psi (s) ds = \sum _{j=0}^N h\rho_j |G(t_i,s_j)|\psi(s_j) +O(h^2), \label{eq:p3.1a}\\ 
\int_0^a |G (\xi_i,s)| \psi (s) ds = \sum _{j=0}^N h\rho_j |G(\xi_i,s_j)|\psi(s_j) +O(h^2), \label{eq:p3.2a}
\end{align}
\end{remark}
\begin{Proposition}\label{prop4}
Under the assumptions of Theorem \ref{thm1} we have the estimates
\begin{align}
\| \Psi _k -\psi_k\|_{\bar{\omega_h}} =O(h^2), \label{est1}\\
\| U _k -u_k\|_{\bar{\omega_h}} =O(h^2), \label{est2},
\end{align}
where $\|. \|_{\bar{\omega_h}}$ is the max-norm of grid function defined on the grid $\bar{\omega_h}$.
\end{Proposition}
\noindent {\it Proof.} We prove the proposition by induction. For $k=0$ we have at once $\| \Psi _0 -\psi_0\|_{\bar{\omega_h}} $ because $\Psi_0(t_i)=f(t_i,0,0)$ and $\psi_0(t_i)=f(t_i,0,0), \; i=\overline{0,N}$, too.
Next, by \eqref{iter2c} and Proposition \ref{prop3} we have
\begin{align*}
u_0(t_i) &= g(t_i)+\int_0^a G(t_i,s)\psi_0(s)ds \\
&=g(t_i)+ \sum_{j=0}^Nh\rho_jG(t_i,s_j)\psi_0(s_j)+O(h^2).
\end{align*}
On the other hand, by \eqref{iter2d} we have
\begin{align*}
U_0(t_i) &=g(t_i)+  \sum _{j=0}^N h\rho_j G(t_i,s_j)\Psi_0(s_j).
\end{align*}
Therefore, 
\begin{align*}
|U_0(t_i)-u_0(t_i)|=O(h^2).
\end{align*}
It implies $\| U_0 -u_0\|_{\bar{\omega_h}}=O(h^2)$. Thus, the estimates \eqref{est1} and \eqref{est2} are valid for $k=0$.\par 
Now, suppose that these estimates are valid for $k\ge 0$. We shall show that they are valid for $k+1$. Indeed, from \eqref{iter3c}, \eqref{iter3d} and the Lipshitz conditions for the function $f(t,u,v)$ we have
\begin{equation}\label{eq:35}
\begin{split}
|\Psi_{k+1}(t_i)-\psi_{k+1}(t_i)|&=|f(t_i, U_k(t_i), V_k(t_i))- f(t_i, u_k(t_i), v_k(t_i)) |\\
& \le L_1 |U_k(t_i)- u_k(t_i)| +L_2 | V_k(t_i)-v_k(t_i) |.
\end{split}
\end{equation}
Now estimate  $| V_k(t_i)-v_k(t_i) |$. We have by Proposition \ref{prop3} 
\begin{align*}
v_k(t_i) &=g(\varphi (t_i))+ \int_0^a G(\varphi (t_i),s)\psi _k(s) ds \\
& =g(\xi_i)+ \sum_{j=0}^N h\rho_j G(\xi_i, s_j)\psi _k (s_j) +O(h^2).
\end{align*}
In view of \eqref{iter2d} we have
\begin{equation}\label{eq:36}
\begin{split}
|V_k(t_i)-v_k(t_i) |&=|\sum_{j=0}^N h\rho_j G(\xi_i, s_j) (\Psi_k(s_j)- \psi_k(s_j))|+O(h^2)\\
& \le \sum_{j=0}^N h\rho_j |G(\xi_i, s_j)| \| \Psi_k-\psi_k\|_{\bar{\omega_h}} +O(h^2).
\end{split}
\end{equation}
Notice that \eqref{eq:p3.2a} for $\psi (s) =1$ gives
\begin{align*}
\int_0^a |G (\xi_i,s)|  ds = \sum _{j=0}^N h\rho_j |G(\xi_i,s_j)| +O(h^2).
\end{align*}
From here it follows that
\begin{equation*}
\begin{split}
&\sum_{j=0}^N h\rho_j |G(\xi_i, s_j)| = \int_0^a |G (\xi_i,s)|  ds +O(h^2)\\
& \le \max_{0\leq t\leq a} \int_0^1 |G(t,s)|ds +O(h^2) =M_0+ O(h^2).
\end{split}
\end{equation*}
Thanks to this estimate, from \eqref{eq:36} we obtain
$$|V_k(t_i)-v_k(t_i) | \le \| \Psi_k-\psi_k\|_{\bar{\omega_h}} + O(h^2).
$$
So, due to the induction hypothesis it implies
\begin{equation}\label{eq:37}
\| V_k-v_k \|_{\bar{\omega_h}}=O(h^2).
\end{equation}
Combining the induction hypothesis $\| U_k-u_k \|_{\bar{\omega_h}}=O(h^2)$ and \eqref{eq:37}, from \eqref{eq:35} we obtain 
\begin{equation}\label{eq:38}
\| \Psi_{k+1}-\psi_{k+1} \|_{\bar{\omega_h}}=O(h^2).
\end{equation}
In order to prove 
\begin{equation}\label{eq:39}
\| U_{k+1}-u_{k+1} \|_{\bar{\omega_h}}=O(h^2).
\end{equation}
we take into account that
\begin{align*}
|U_{k+1}(t_i)-u_{k+1}(t_i)| \le \sum_{j=0}^Nh\rho_j |G(t_i,s_j)|  |\Psi_{k+1}(s_j)-\psi_{k+1}(s_j) | +O(h^2).
\end{align*}
Doing the similar argument as above and using the proved estimate \eqref{eq:38} it is easy to obtain
\begin{align*}
|U_{k+1}(t_i)-u_{k+1}(t_i)| = O(h^2),
\end{align*}
or \eqref{eq:39}.\\
Thus, Proposition \ref{prop4} is proved.\par 
Now combining  Proposition \ref{prop4} with Theorem \ref{thm2} we obtain the following result.
\begin{Theorem}\label{thm3}
Under the assumptions of Theorem \ref{thm1} for the approximate solution of the problem \eqref{eq:2}-\eqref{eq:3} obtained by the discrete iterative method \eqref{iter1d}-\eqref{iter3d} we have the estimate
\begin{equation*}
\|U_k-u\|_{\bar{\omega_h}} \leq M_0p_kd +O(h^2), 
\end{equation*}
where $p_k$ and $d$ are defined by \eqref{eqpkd}.
\end{Theorem}

\begin{remark}
The results in Section \ref{sec2} and \ref{sec3} are obtained for the nonlinear third order FDE with nonlinear term $f=f(t,u(t),u(\varphi(t)))$. Analogously, it is possible to obtain similar results of existence and convergence of the iterative method at continuous level for the general case $$f=f(t,u(t),u(\varphi(t)), u'(\varphi_1(t)), u''(\varphi_2(t))).$$  But for numerical realization of the iterative method it is needed to take attention that the second derivative $\frac{\partial ^2 G(t,s)}{\partial t^2}$ of the Green function has discontinuity at the line $s=t$. In this case for computing integrals containing $\frac{\partial  G(t,s)}{\partial t}$ and $\frac{\partial ^2 G(t,s)}{\partial t^2}$ it is needed to use the formulas constructed in our previous work \cite{Dang2}.
\end{remark}
\begin{remark}
The technique developed in this paper for the third order nonlinear FDE can be applied to nonlinear FDE of any order.
\end{remark}

\section{Examples}
In all numerical examples below we perform the iterative method \eqref{iter1d}-\eqref{iter3d} until $\| \Psi _k -\Psi_{k-1}\|_{\bar{\omega_h}} \le 10^{-10}$. In the tables of results for the convergence of the iterative method $Error=\|U_K-u \|_{\bar{\omega_h}}$, $K$ is the number of iterations performed.
 
\noindent {\bf Example 1.} Consider the following problem
\begin{equation}\label{eq:exam1}
\begin{aligned}
u'''(t)&=e^t-\frac{1}{4}u(t)+\frac{1}{4} u^2(\frac{t}{2}), \quad 0<t<1,\\
u(0)&=1,\; u'(0)=1, \; u'(1)=e
\end{aligned}
\end{equation}
with the exact solution $u(t)=e^t.$
The Green function for the above problem is
\begin{equation*}
\begin{aligned}
G(t,s)=\left\{\begin{array}{ll}
\dfrac{s}{2}(t^2-2t+s), \quad 0\le s \le t \le 1,\\
\, \, \dfrac{t^2}{2}(s-1), \quad 0\le t \le s \le 1.\\
\end{array}\right.
\end{aligned}
\end{equation*}
So, we have 
\begin{equation*}
M_0  =\max_{0\leq t\leq a} \int_0^1 |G(t,s)|ds =\frac{1}{2} .
\end{equation*}
The second degree polynomial satisfying the boundary conditions of the problem is
$$g(t)=1+t+\frac{e-1}{2} t^2.$$
Therefore, $\|g\|=2+\dfrac{e-1}{2}=2.7183$. In this example $f(t,u,v)=e^t-\frac{1}{4}u+\frac{1}{4}v^2  $.  It is easy to verify that for $M=6.5$ we have $|f(t,u,v)|\le M$ in the domain $\mathcal{D}_M$ defined by \eqref{eq:11}. Moreover, in this domain the function $f(t,u,v)$ satisfies the Lipshitz conditions in $u$ and $v$ with the coefficients $L_1=\frac{1}{4}$ and $L_2=1.7004$. Therefore, $q:=(L_1+L_2)M_0= 0.16$. Thus, all the assumptions of Theorem \ref{thm1} are satisfied. By the theorem, the problem \eqref{eq:exam1} has a unique solution. This is the above exact solution.\par 
The results of convergence of the iterative method \eqref{iter1d}-\eqref{iter3d} are given in Table \ref{table1}.
\begin{table}[ht!]
\centering
\caption{The convergence in Example 1. }
\label{table1}
\begin{tabular}{cccc}
\hline 
$N$ &	$h^2$ &	$K$	&$Error$ \\
\hline 
50	& 4.0000e-04 &	3 &	6.1899e-05\\
100	& 1.0000e-04 &	3&	1.5475e-05\\
150	&4.4444e-05	&   3&6.877 -06\\
200	&2.5000e-05	&3	&3.8688e-06\\
300	&1.1111e-05	&3	&1.7195e-06\\
400&	6.2500e-06&	3	&9.6721e-07\\
500	&4.0000e-06	&3	&6.1901e-07\\
800	&1.5625e-06&	3&	6.1901e-07\\
1000&	1.0000e-06	&3	&1.5475e-07
\\
\hline 
\end{tabular} 
\end{table}
From  this table we see that the results of computation support the conclusion that the accuracy of the iterative method is $O(h^2)$.
\begin{remark}
Theorem \ref{thm3} gives sufficient conditions for convergence of the iterative method \eqref{iter1d}-\eqref{iter3d}. In the cases when these conditions are not satisfied the iterative also can converge to some solution. For example, for the case $f(t,u,v)=e^t+u^2+v^2+1$ with the same boundary conditions as in \eqref{eq:exam1} the iterative method converges after 15 iterations. And for the case $f(t,u,v)=e^{2t}-u^3+v^2+5$ after $16$ iterations the iterative process reaches the $TOL=10^{-10}$. Notice that the number of iterations do not depend on the grid size as in Example 1.
\end{remark}

\begin{figure}[ht]
\begin{center}
\includegraphics[height=6cm,width=9cm]{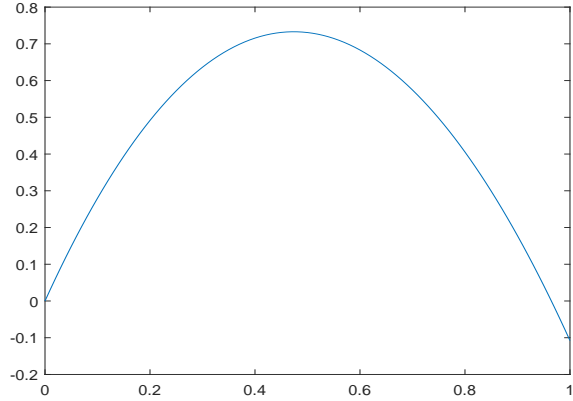}
\caption{The graph of the approximate solution in Example $2$. }
\label{fig1}
\end{center}
\end{figure}

\noindent {\bf Example 2.} Consider the following problem
\begin{equation}\label{eq:exam2}
\begin{aligned}
u'''(t)&=\sin (u^2(t)) +\cos (u^2 (t^2)), \quad 0<t<1,\\
u(0)&=0,\; u'(0)=\pi, \; u'(1)=-\pi .
\end{aligned}
\end{equation}
For this problem $f(t,u,v)=\sin (u^2)+\cos (v^2)$ and $\varphi (t)=t^2$. It is easy to verify that all the conditions of Theorem \ref{thm2} are satisfied, therefore the problem has a unique solution. Also, by Theorem \ref{thm3} the iterative method \eqref{iter1d}-\eqref{iter3d} converges. With $TOL=10^{-10}$ the iterative method for any number of grid points stops after $8$ iterations. The graph of the approximate solution is depicted in Figure \ref{fig1}.

\noindent {\bf Example 3.} (Example 5 in \cite{Khuri}) Consider the problem
\begin{equation}\label{eq:exam3}
\begin{aligned}
u'''(t)&=-1+2u^2(t/2), \quad 0<t<\pi,\\
u(0)&=0,\; u'(0)=1, \; u(\pi)=0.
\end{aligned}
\end{equation}
which has exact solution $u(t)=\sin (t)$. For this problem the Green function has the form
\begin{equation*}
\begin{aligned}
G(t,s)=\left\{\begin{array}{ll}
-\dfrac{t^2 (\pi -s)^2}{2\pi ^2}+ \dfrac{(t-s)^2}{2}, \quad 0\le s \le t \le \pi,\\
\, \, -\dfrac{t^2 (\pi -s)^2}{2\pi ^2}, \quad 0\le t \le s \le \pi\\
\end{array}\right.
\end{aligned}
\end{equation*}
and the function $f(t,u,v)=-1+2v^2$.

The results of convergence of the iterative method \eqref{iter1d}-\eqref{iter3d} for this example are given in Table \ref{table2}
\begin{table}[ht!]
\centering
\caption{The convergence in Example 3. }
\label{table2}
\begin{tabular}{cccc}
\hline 
$N$ &	$h^2$ &	$K$	&$Error$ \\
\hline 
50	& 4.0000e-04 &	25 &	1.4455e-04\\
100	& 1.0000e-04 &	25&	3.6142e-05\\
150	&4.4444e-05	&   25& 1.6063e-05\\
200	&2.5000e-05	&25	& 9.0345e-06\\
300	&1.1111e-05	&25	&4.0155e-06\\
400&	6.2500e-06&	25	&2.2587e-06 \\
500	&4.0000e-06	&25	&1.4456e-06\\
800	&1.5625e-06& 25&	5.6467e-07\\
1000&	1.0000e-06	&25	&3.6139e-07\\
\hline 
\end{tabular} 
\end{table}
From Table \ref{table2} it is seen also that the numerical method has the accuracy $O(h^2)$.

\section{Conclusion}
In this paper we have proposed a unified approach to nonlinear functional differential equations via boundary value problems for nonlinear third order functional differential equations as a particular case. We have established the existence and uniqueness of solution and proved the convergence of the discrete iterative method for finding the solution. Some examples demonstrate the validity of the theoretical results and the efficiency of the numerical method.\par 
The proposed approach can be applied to boundary value problems for nonlinear functional differential equations of any order. It also can be applied to integro-differential equations.

\end{document}